\documentclass[11pt]{amsart}
\usepackage{epsfig}
\usepackage{graphics}

\newtheorem{theorem}{Theorem}
\newtheorem{lemma}[theorem]{Lemma}

\newtheorem{corollary}[theorem]{Corollary}

\theoremstyle{definition}

\theoremstyle{remark}

 \def\b{{\mathcal{B}}}
 
 \def\Z{{\mathbb{Z}}}

 \def\S{{\Sigma}}
\def\mod{{\rm Mod}}
\def\emod{{\rm Mod}^*}

\begin{document}

\newenvironment{prooff}{\medskip \par \noindent {\it Proof}\ }{\hfill
$\square$ \medskip \par}
    \def\sqr#1#2{{\vcenter{\hrule height.#2pt
        \hbox{\vrule width.#2pt height#1pt \kern#1pt
            \vrule width.#2pt}\hrule height.#2pt}}}
    \def\square{\mathchoice\sqr67\sqr67\sqr{2.1}6\sqr{1.5}6}
\def\pf#1{\medskip \par \noindent {\it #1.}\ }
\def\endpf{\hfill $\square$ \medskip \par}
\def\demo#1{\medskip \par \noindent {\it #1.}\ }
\def\enddemo{\medskip \par}
\def\qed{~\hfill$\square$}

 \title[Generating the surface mapping class group by two elements]
 {Generating the surface mapping class group by two elements}

\author[Mustafa Korkmaz]{Mustafa Korkmaz $\,^*$}

 \address{Department of Mathematics, Middle East Technical University,
 06531 Ankara, Turkey} \email{korkmaz@arf.math.metu.edu.tr}

 \date{\today}
 \thanks{* Supported by T\"UBA-GEB\.IP}
\begin{abstract}
Wajnryb proved in~\cite{w1} that the mapping class group of an orientable surface
is generated by two elements. We prove that one of these generators can be taken
as a Dehn twist. We also prove that the extended mapping class group is generated
by two elements, again one of which is a Dehn twist. Another result we prove is that the 
mapping class groups are also generated by two elements of finite order.
\end{abstract}

 \maketitle
  \setcounter{secnumdepth}{1}
 \setcounter{section}{0}

\section{Introduction}
Let $\S$ be a compact connected oriented surface of genus $g$ with
one boundary component. We denote by $\mod_g^1$ the mapping class
group of $\S$, the group of isotopy classes of
orientation-preserving diffeomorphisms $\S\to \S$ which restrict to
the identity on the boundary. The isotopies are also required to
fix the points on the boundary. If the diffeomorphisms and the
isotopies are allowed to permute the points on the boundary of $\S$, 
then we get the group $\mod_{g,1}$. The extended mapping
class group $\emod_{g,1}$ of $\S$ is defined to be the group of
isotopy classes of all (including orientation-reversing)
diffeomorphisms of $\S$. These three groups are related to each
other as follows: $\mod_{g,1}$ is contained in $\emod_{g,1}$ as a
subgroup of index two and the groups $\mod_g^1$ and $\mod_{g,1}$
fit into a short exact sequence
$$1\rightarrow \Z \rightarrow  \mod_g^1 \rightarrow \mod_{g,1} \rightarrow 1,$$
where $\Z$ is the subgroup of  $\mod_g^1$
generated by the Dehn twist about a simple closed curve
parallel to the boundary component of $\S$.

In this paper, we will be interested in the groups $\mod_g^1$, $\mod_{g,1}$ and $\emod_{g,1}$.
The mapping class group and the extended mapping class group
of the closed surface of genus $g$ obtained from $\S$ by gluing a disc
along the boundary component are denoted by $\mod_g$ and $\emod_g$.

\begin{figure}[hbt]
 \begin{center}
    \includegraphics[width=10cm]{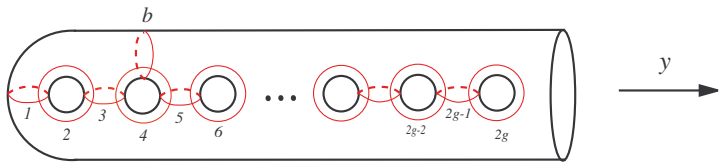}
  \caption{The curve labelled $i$ is $a_i$.}
  \label{figure1}
   \end{center}
 \end{figure}

The mapping class group is a central object in low-dimensional
topology. Therefore, its algebraic structures are of interest. The
problem of finding the generators for the mapping class group of a
closed orientable surface was first considered by Dehn. He proved
in~\cite{d} that $\mod_g$ is generated by a finite set of Dehn
twists. About quarter century later, Lickorish~\cite{l1,l2} also
proved the same result, showing that $3g-1$ Dehn twists generate
$\mod_g$. This number was improved to $2g+1$ by
Humphries~\cite{h}. These $2g+1$ generators are the Dehn twists
about the curves $b,a_1,a_2,\ldots, a_{2g}$ in
Figure~\ref{figure1}, where the closed surface is obtained from
$\S$ by gluing a disc along the boundary component. Humphries
proved, moreover, that in fact the number $2g+1$ is minimal; i.e.
$\mod_g$ cannot be generated by $2g$ (or less) Dehn twists.
Johnson~\cite{j} proved that the $2g+1$ Dehn twists about
$b,a_1,a_2,\ldots, a_{2g}$ on $\S$ also generate $\mod_g^1$.
Finally, the minimal number of generators for the mapping class
group is determined by Wajnryb~\cite{w1}. He showed
 that $\mod_g^1$, and hence $\mod_g$,  can be generated
by two elements; one is a product of two Dehn twists (one is right
and one is left) and the other is a product of $2g$ Dehn twists.
Since the mapping class group is not abelian, the number two is
minimal. Recently, it was shown by Brendle and Farb in~\cite{bf}
that the mapping class group is generated by three torsion
elements and by seven involutions.

Since $\mod_{g,1}$ is a quotient of $\mod_g^1$, it is generated by the corresponding
$2g+1$ Dehn twists. In order to generate the extended mapping class group $\emod_{g,1}$,
it suffices to add one more generator, namely the isotopy class of any orientation-reversing diffeomorphism.

In this paper we have three main results. First, we improve Wajnryb's result.
We show that one of the two generators of $\mod_g^1$ can be taken as a Dehn twist.
All Dehn twists involved in our generators are Dehn-Lickorish-Humphries generators.
We also prove that the extended mapping class group $\emod_{g,1}$ is
generated by two elements, again one of which is a Dehn twist.
Our proof is independent from that of Wajnryb~\cite{w1}. 
 Next, we prove that the mapping class groups $\mod_{g,1}$ and $\mod_g$
are also generated by two torsion elements. Of course, this is not  true
for $\mod_g^1$ since it is torsion-free. In the last section of the paper,
we transform the presentation of the mapping class group in~\cite{w2} into a presentation
on our two generators.

\section{Preliminaries}

Recall that if $a$ is a simple closed curve on the oriented surface
$\S$, then the (right) Dehn twist
$A$ about $a$ is the isotopy class of the diffeomorphism obtained by cutting $\S$
along $a$, twisting one of the side by $2\pi$ to the right and gluing two sides of $a$ back to each other.
We denote the curves by the letters $a,b,c,d$ (possibly with subscripts) and the Dehn twists
about them by the corresponding capital letters $A,B,C,D$.
Notationally we do not distinguish a diffeomorphism/curve and its isotopy class.

We use the functional notation for the composition of two
diffeomorphisms; if $F$ and $G$ are two diffeomorphisms, then the
composition $FG$ means that $G$ is applied first.

\begin{figure}[hbt]
 \begin{center}
    \includegraphics[width=8cm]{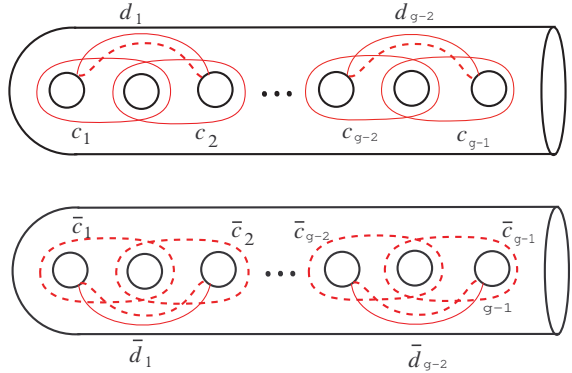}
  \caption{The curves $c_i,d_i,\bar{c}_i$ and $\bar{d}_i$.}
  \label{figure2}
   \end{center}
 \end{figure}

We define the curves $c_i$, $d_i$, $\bar{c}_i$ and
$\bar{d}_i$ on $\S$ as shown in Figure~\ref{figure2}, so that
$\bar{c}_i$ (resp. $\bar{d}_i$) is obtained from $c_i$
(resp. $d_i$) by rotating the surface $\S$ about the $y$-axis by
$\pi$. These curves will be used through out the paper.
Here, we assume that the surface of the paper is the $yz$-plane,
the positive side of the $x$-axis is pointing above the page and
the surface is invariant under the rotation by $\pi$ about the $y$-axis.

Let $G$ be a subgroup of $\mod_g^1$ or $\emod_{g,1}$.
Then $G$ acts on the set of the isotopy classes of simple closed curves.
If $c$ is a simple closed curve, then we denote by
$c^G$ the $G-$orbit of $c$;
$$c^G=\{ \, F(c) : F\in G \}. $$

We record the following lemmas.

\begin{lemma} \label{lemmaxx}
Let $c$ be a simple closed curve on $\Sigma$, let $F$ be a self-diffeomorphism  of $\Sigma$ 
and let $F(c)=d$. Then  $FCF^{-1}=D^r$, where $r=\pm 1$ depending on whether $F$ is 
orientation-preserving or orientation-reverving.
 \end{lemma}

\begin{lemma} 
Let $c$ and $d$ be two simple closed curves on $\Sigma$. If $c$ is disjoint from $d$,
then $CD=DC$. 
 \end{lemma}

\begin{lemma} \label{lemma1}
Let $c$ and $d$ be two simple closed curves on $\Sigma$.
Suppose that $C\in G$ and the curve $d$ is contained in
$c^G$. Then $D$ is also contained in $G$.
 \end{lemma}
\begin{proof}
Since $d\in c^G$, $F(c)=d$ for some $F\in G$. By Lemma~\ref{lemmaxx}, 
$F C F^{-1} =D^r$. This shows that $D$ is contained in $G$.
\end{proof}

\section{The mapping class group $\mod_g^1$}
\label{sec:mcg}

Let $S$ denote the product $A_{2g}A_{2g-1}\cdots A_{2}A_{1}$ of $2g$ Dehn twists
in $\mod_g^1$ and let $G$ be the subgroup of $\mod_g^1$ generated by $B$ and $S$. We
prove in this section that $G=\mod_g^1$. It follows that the
mapping class groups $\mod_{g,1}$ and $\mod_g$ are also generated
by $B$ and $S$. The main idea of the proof is to show that the
$G$-orbit $b^G$ of the curve $b$ contains the simple closed curves $a_1,a_2,\ldots,a_{2g}$.

\begin{lemma} \label{lemma2}
The curves $c_1,c_2,\ldots ,c_{g-1}$ and $d_1,d_2,\ldots ,d_{g-2}$
of Figure~$\ref{figure2}$ are contained in $b^G$, the $G$-orbit of $b$.
 \end{lemma}
\begin{proof}
This follows from these easily verified facts:
The diffeomorphism $S^{-1}$ maps

$\bullet$ $b$ to $c_1$,

$\bullet$  $c_i$ to $d_i$, and

$\bullet$ $d_i$ to $c_{i+1}$.
\end{proof}

\noindent
{\bf Remark.} It can be shown that $S^{2g+1}(c_i)=\bar{c}_i$ and 
$S^{2g+1}(d_i)=\bar{d}_i$, so that $\bar{c}_i$ and $\bar{d}_i$ are also in $b^G$.
These facts will be used in Section~\ref{sec:torsion} in the proof of the fact that the mapping class groups
$\mod_{g,1}$ and $\mod_{g}$ are generated by two torsion elements

\begin{figure}[hbt]
 \begin{center}
    \includegraphics[width=13cm]{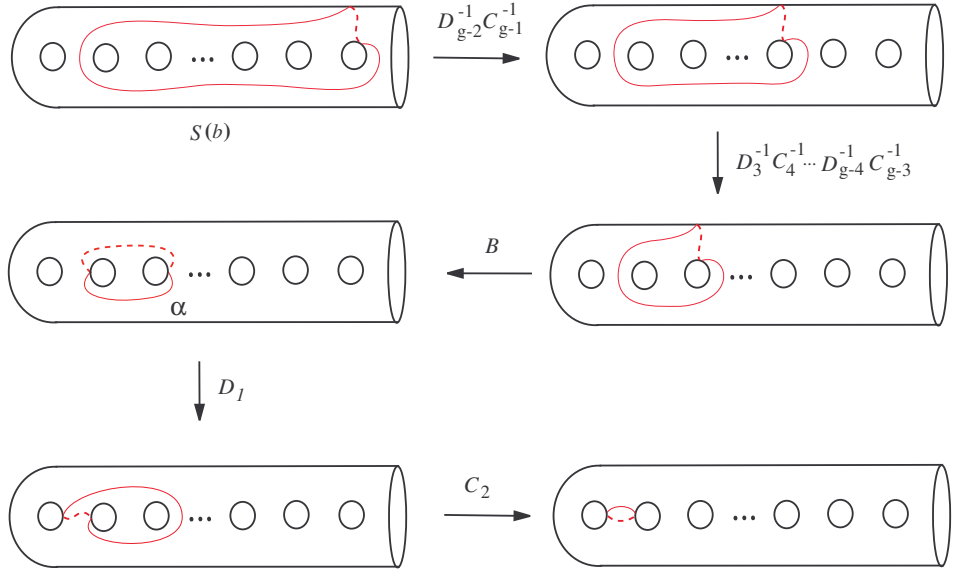}
  \caption{$g$ is odd.}
  \label{alpha}
  \end{center}
 \end{figure}

 \begin{theorem} \label{thmmcg}
 Suppose that $g\geq 2$. The subgroup $G$ generated by $B$ and $S$
 is equal to the mapping class group $\mod_g^1$.
 \end{theorem}
\begin{proof}
It can easily be shown that $S(a_i)=a_{i-1}$. Hence
$SA_iS^{-1}=A_{i-1}$, and thus $A_{i-1}\in G$ if and only if
$A_i\in G$.

Suppose that the curve $a_{i_0}$ is contained in $b^G$ for some
$i_0$. Since $B\in G$, by Lemma~\ref{lemma1} we get
that $A_{i_0}$ is also contained in $G$. Therefore, all $A_i$ are
contained in $G$. Since the mapping class group $\mod_g^1$ is
generated by the $2g+1$ Dehn twists $A_1,A_2,\ldots, A_{2g}$ and
$B$, we conclude that the subgroup $G$ is in fact equal to
$\mod_g^1$. Therefore, in order to finish the proof of the theorem
it suffices to show that $a_i$ is contained in $b^G$ for some $i$
with $1\leq i\leq 2g$.

Suppose first that $g$ is odd. It is easy to see that the
diffeomorphism
$$BD_3^{-1}C_4^{-1}D_5^{-1}C_6^{-1}\cdots D_{g-2}^{-1}C_{g-1}^{-1}$$
maps the curve $S(b)$ to a curve $\alpha$ and $C_2D_1$ maps
$\alpha$ to $a_3$ (cf. Figure~\ref{alpha}). Since all $C_i$ and $D_i$ are in $G$,
the curve $a_3$ is contained in the $G$-orbit $b^G$ of $b$.

Suppose now that $g$ is even. In this case the diffeomorphism
$$B D_2^{-1}C_3^{-1}D_4^{-1}C_5^{-1}\cdots D_{g-2}^{-1}C_{g-1}^{-1}$$
maps the curve $S(b)$ to $a_4$. Again, since all $C_i$ and $D_i$
are in $G$, we conclude that $a_4$ is contained in $b^G$.

This concludes the proof of the theorem.
\end{proof}

\section{The extended mapping class group $\emod_{g,1}$}

In this section we prove that the extended mapping class group $\emod_{g,1}$ is
also generated by two elements, one of which is a Dehn twist.

Consider the surface $\S$ embedded in the $3$-space as shown in
Figure~\ref{figure1}. Let $R$ denote the reflection across the
$xy$-plane and let $T$ denote the product 
$A_{2g}A_{2g-1}\cdots A_{2}A_{1}R$. 
Let $H$ denote the subgroup of $\emod_{g,1}$
generated by $B$ and $T$. We prove that $H=\emod_{g,1}$. Again, it
follows that $\emod_{g}$ is also generated by $B$ and $T$. Recall
that the $H$-orbit of the simple closed curve $b$ is denoted by
$b^H$.

 \begin{lemma} \label{lemmaxx}
{\rm (i)} If $g$ is even, then the curves
$c_1,\bar{c}_2,c_3,\bar{c}_4,\ldots ,\bar{c}_{g-2},c_{g-1}$
and $\bar{d}_1,d_2,\bar{d}_3,d_4,\ldots, \bar{d}_{g-3},d_{g-2}$
of Figure~$\ref{figure2}$ are contained in $b^H$.

{\rm (ii)} If $g$ is odd, then the curves
$c_1,\bar{c}_2,c_3,\bar{c}_4,\ldots ,c_{g-2},\bar{c}_{g-1}$
and $\bar{d}_1,d_2,\bar{d}_3,$ $d_4,\ldots ,d_{g-3},\bar{d}_{g-2}$
of Figure~$\ref{figure2}$ are contained in $b^H$.
 \end{lemma}
\begin{proof}
It can be shown easily that the diffeomorphism $T^{-1}$ maps

$\bullet$ $b$ to $c_1$,

$\bullet$ $c_i$ to $\bar{d}_i$,

$\bullet$ $\bar{c}_i$ to $d_i$,

$\bullet$ $d_i$ to $c_{i+1}$, and

$\bullet$ $\bar{d}_i$ to $\bar{c}_{i+1}$.\\

The lemma follows from these and Lemma~\ref{lemma1}.
\end{proof}

\begin{figure}[hbt]
 \begin{center}
    \includegraphics[width=13cm]{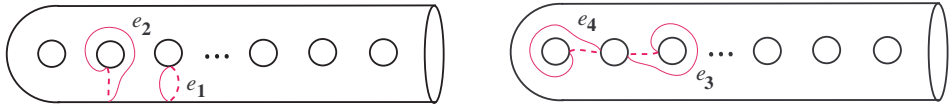}
  \caption{The curves $e_i$.}
  \label{e_i}
   \end{center}
 \end{figure}

\begin{figure}[hbt]
 \begin{center}
    \includegraphics[width=13cm]{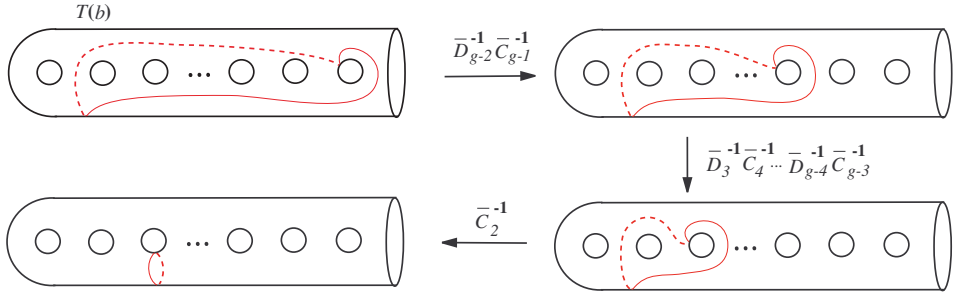}
  \caption{The proof of $e_1\in b^H$.}
  \label{e_1}
   \end{center}
 \end{figure}

\begin{figure}[hbt]
 \begin{center}
    \includegraphics[width=13cm]{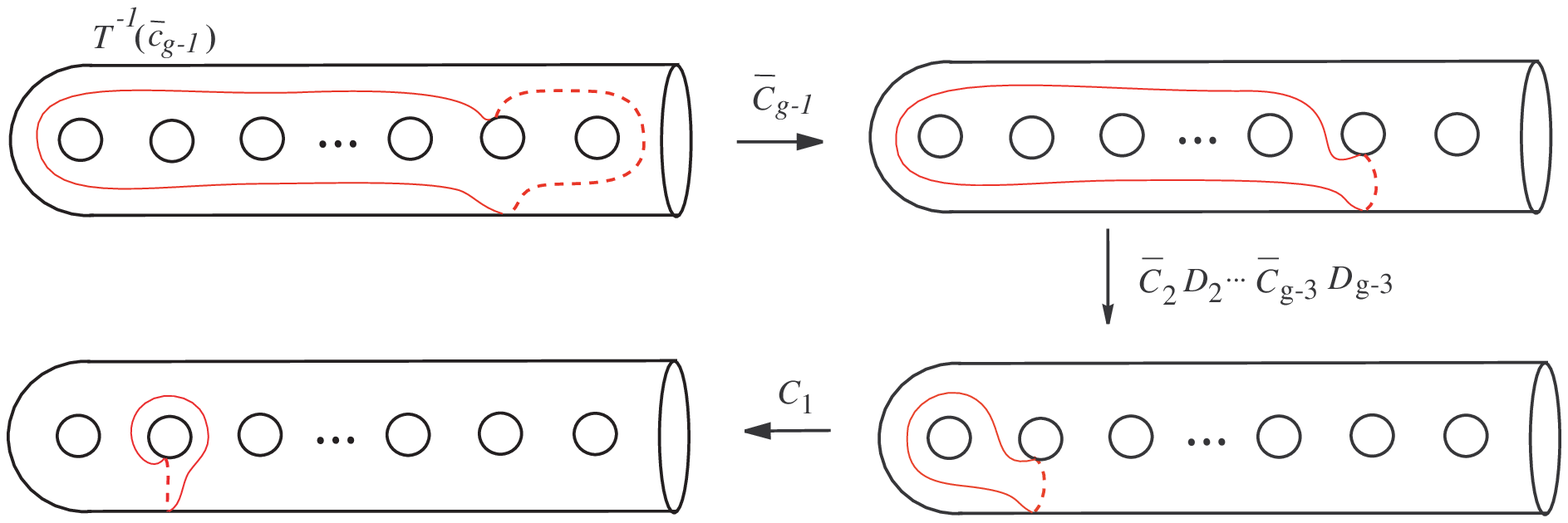}
  \caption{The proof of $e_2\in b^H$.}
  \label{e_2}
   \end{center}
 \end{figure}

\begin{figure}[hbt]
 \begin{center}
    \includegraphics[width=13cm]{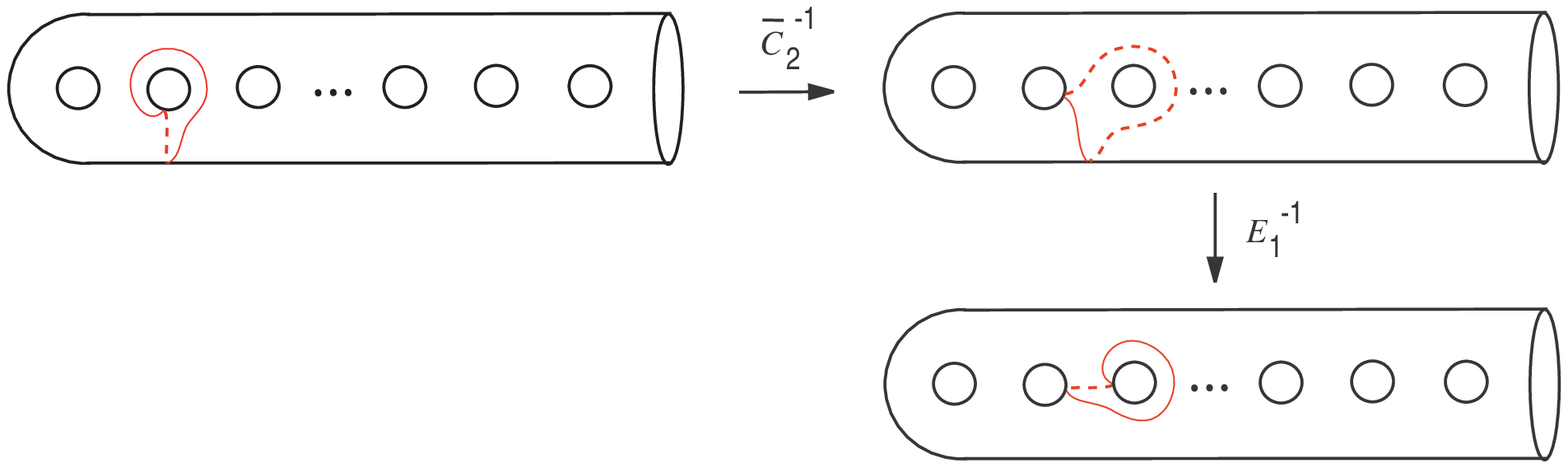}
  \caption{The proof of $e_3\in b^H$.}
  \label{e_3}
   \end{center}
 \end{figure}

\begin{figure}[hbt]
 \begin{center}
    \includegraphics[width=13cm]{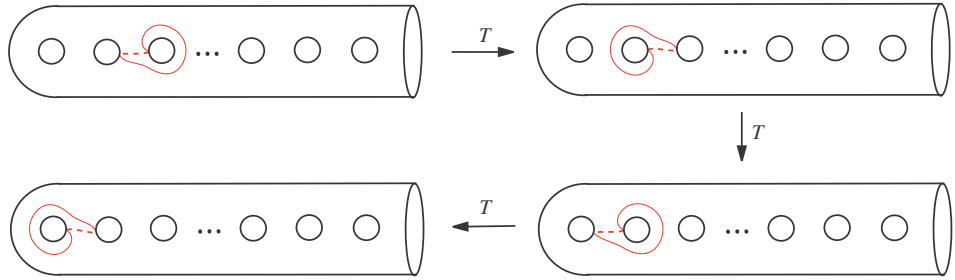}
  \caption{The proof of $e_4\in b^H$.}
  \label{e_4}
   \end{center}
 \end{figure}

 \begin{lemma} \label{lemma4}
Suppose that $g$ is odd. Then the curves
$e_1,e_2,e_3$ and $e_4$ in Figure~$\ref{e_i}$ are contained in $b^H$.
 \end{lemma}
\begin{proof}
Let $U_i$ denote $(\bar{C}_{i})^{-1}(\bar{D}_{i+1})^{-1}$.
If $i$ is even then $U_i$ is contained in $H$ by Lemma~\ref{lemmaxx}.
Now the diffeomorphism
$$ U_2U_4U_6\cdots U_{g-3}(\bar{C}_{g-1})^{-1}T$$
is contained in $H$ and it maps
the curve $b$ to $e_1$ (c.f. Figure~\ref{e_1}). This proves that $e_1\in b^H$.

For the proof of $e_2\in b^H$ let $U$ denote the diffeomorphism
$$ \bar{C}_{2} D_2 \bar{C}_{4} D_4\cdots \bar{C}_{g-3}
D_{g-3}.$$
Now it suffices to show that
$$ C_1U\bar{C}_{g-1}T^{-1} $$
maps the curve $\bar{c}_{g-1}$ to $e_2$. This can be seen
from Figure~\ref{e_2}.

The curve $e_3$ is the image of the curve $e_2$ under the
diffeomorphism $E_1^{-1}(\bar{C}_{2})^{-1}$ (c.f.
Figure~\ref{e_3}). This proves that $e_3\in b^H$.

Finally, $T^3(e_3)=e_4$ proving that $e_4\in b^H$ (c.f.
Figure~\ref{e_4}).
\end{proof}

\begin{figure}[hbt]
 \begin{center}
    \includegraphics[width=13cm]{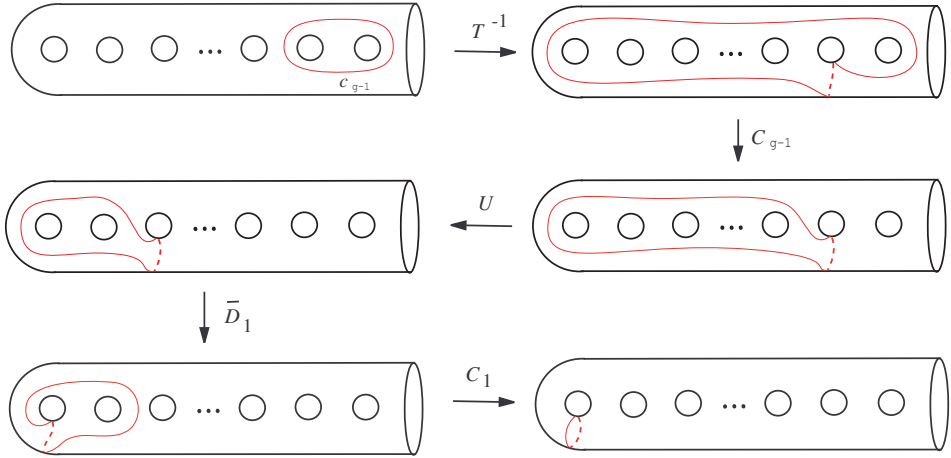}
  \caption{The proof of $a_1\in b^H$ for even $g$.}
  \label{extmcgeven}
   \end{center}
 \end{figure}

 \begin{theorem} \label{even-g}
 The subgroup $H$ generated by $B$ and $T$
 is equal to the extended mapping class group $\emod_{g,1}$.
 \end{theorem}
\begin{proof}
We prove this theorem in the same way as Theorem~\ref{thmmcg}; we
show that the $H$-orbit $b^H$ of the curve $b$ contains simple closed curves
$a_1,a_2,\ldots, a_{2g}$.

 It is easy to show that $T(a_{i})=a_{i-1}$ for
$i=2,3,\ldots, 2g-1$. Hence, $TA_iT^{-1}=A_{i-1}$, and thus $A_i\in
H$ if and only if $A_{i-1}\in H$.

Suppose that $a_{i_0}\in b^H$ for some
$i_0$. Since $B$ is contained in $H$, by Lemma~\ref{lemma1} we get
that $A_{i_0}$ is also contained in $H$. Therefore, all $A_i$ are
contained in $H$. Since $T\in H$, the reflection $R$ is also
contained in $H$. The extended mapping class group
$\emod_{g,1}$ is generated by the $2g+1$ Dehn twists
$B,A_1,A_2,\ldots, A_{2g}$ and the reflection $R$. We conclude
that the subgroup $H$ is in fact equal to $\emod_{g,1}$.
Therefore, in order to finish the proof of the theorem it suffices
to show that $b^H$ contains $a_i$ for some $i$ with $1\leq i\leq
2g$.

Suppose first that $g$ is even.
It follows from Lemma~\ref{lemma1} and Lemma~\ref{lemmaxx} that
the diffeomorphism
$$V=C_1\bar{D}_{1}C_3\bar{D}_{3}\cdots
C_{g-3}\bar{D}_{g-3}C_{g-1}T^{-1}$$ is contained in $H$.
Figure~\ref{extmcgeven} shows that $V$ maps the curve $c_{g-1}$ to
$a_1$. Since $c_{g-1}\in b^H$, $a_1$ is also in $b^H$. This
finishes the proof in this case.

Suppose now that $g$ is odd. It easy to verify that
$(C_1)^{-1}E_4(e_2)=a_1$. Since $e_2\in b^H$, and $C_1$ and $E_4$
are in $H$, we conclude that $a_1\in b^H$.

This completes the proof of the theorem.
\end{proof}

 \section{Generating mapping class groups by two elements of finite order}
\label{sec:torsion}

In this section we prove that the mapping class groups $\mod_{g,1}$ and $\mod_{g}$ are 
generated by two elements of finite order. Clearly, this will be the minimum number of
such generators. In their paper~\cite{bf}, Brendle and Farb 
proved that these mapping class groups are generated by three torsion elements and asked if 
they can be generated by two. 
Therefore our result gives a positive answer to their question.

Let $\Sigma$  be a surface with one boundary component as in Figure~\ref{figure1}. 
In $\mod_{g,1}$, let $S$ denote the product $A_{2g}A_{2g-1}\cdots A_2A_1$. Note that
$S$ is of order $4g+2$. Throughout this section, let $G$ denote the subgroup of $\mod_{g,1}$
generated by the two torsion elements $S$ and $BSB^{-1}$. We prove that $G=\mod_{g,1}$
for $g\geq 3$. 

\begin{figure}[hbt]
 \begin{center}
    \includegraphics[width=8cm]{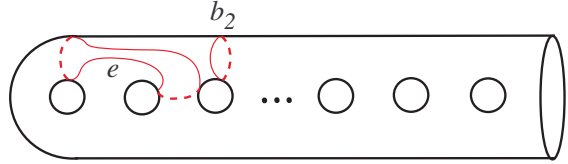}
  \caption{The curves $e$ and $b_2$.}
  \label{lantern}
   \end{center}
 \end{figure}

In the proof of the main result of this section, we use the celebrated lantern relation,
which was discovered by Dehn and redisovered by Johnson. This relation is read as follows;
    \begin{eqnarray} \label{lant1}
    A_1A_3A_5B_2=BD_1E,
   \end{eqnarray}
where the $B_2$ and $E$ are the Dehn twists about the curves $b_2$ and $e$ in Figure~\ref{lantern}.
We rewrite the relation as
   \begin{eqnarray} \label{lant2}
   A_1=(BA_3^{-1}) ( D_1A_5^{-1} ) (EB_2^{-1}).
   \end{eqnarray}
The strategy of the proof is to show that the statements inside each paranthesis
are in $G$. Then the rest of the proof is easy.

Let us define $\b$ as the subset of nonseparating simple closed curves consisting
of those curves $x$ such that $XB^{-1}\in G$. That is,
$$
\b=\{ \ x \ | \ \mbox{$x$ is a nonseparating simple closed curve and} \ XB^{-1}\in G     \}.   
$$

We first state the following easy to prove, but very useful, lemmas.

\begin{lemma}
Suppose that two simple closed curves $x$ and $y$ are contained in $\b$. Then $XY^{-1}$ is contained
in the subgroup $G$.
\end{lemma}

\begin{lemma}
Suppose that $y\in \b$ and $XY^{-1}\in G$. Then $x\in \b$.
\end{lemma}

\begin{lemma}\label{lemma9}
Suppose that a curve $y$ is contained in $\b$ and $x$ is in the $\langle S \rangle$-orbit of $y$.
Then $x$ is contained in $\b$. 
\end{lemma}
\begin{proof}
By assumption, there is an integer $k$ such that $x=S^k(y)$. 
Since the element 
$$
XY^{-1}
=(S^kYS^{-k})Y^{-1}
=S^k (YB^{-1}) (BSB^{-1})^{-k} (BY^{-1})
$$ 
is contained in $G$,
the lemma is proved.
\end{proof}

\begin{corollary}
The curves $c_i,\bar{c}_i,d_j,\bar{d}_j, S(b)$ are contained in $\b$ for all 
$i=1,2,\ldots, g-1$ and $j=1,2,\ldots, g-2$. 
\end{corollary}
\begin{proof}
The fact that $c_i,d_j$ are in the $\langle S \rangle$-orbit of $b$ is shown in Section~\ref{sec:mcg}.
It can be shown that $S^{2g+1}$ is isotopic to the rotation by $\pi$ about the $y$-axis. Since
 $\bar{c}_i=S^{2g+1}(c_i )$ and $\bar{d}_j=S^{2g+1}(d_j )$, the corollary follows 
from Lemma~\ref{lemma9}.
\end{proof}

\begin{lemma} \label{tor}
Suppose that $g\geq 3$. Then for each $i=1,2,\ldots, g$, the curves $a_i$, $b_2$ and $e$
are contained in $\b$. 
\end{lemma}
\begin{proof}
We prove first that the curves $a_i$ are contained in $\b$. Since all $a_i$ are contained in the same 
$\langle S\rangle$-orbit, it is enough to prove that any one of $a_i$ is in $\b$.

Suppose that $g$ is even, so that $g\geq 4$. Let $V$ denote the product 
$C_3^{-1}D_4^{-1}\cdots D_{g-4}^{-1}C_{g-3}^{-1}D_{g-2}^{-1}C_{g-1}^{-1}$. We have shown 
in the proof of Theorem~\ref{thmmcg} that the diffeomorphism $BD_2^{-1}VS$ maps the curve $b$ 
to $a_4$. Since the curve $\bar{c}_1$ is disjoint from each $c_i$ and $d_j$ for $i\geq 1$, $j\geq 3$,
the Dehn twist $\bar{C}_1$ about $\bar{c}_1$ commutes with each $C_i$ and $D_j$. 
Therefore it commutes with $V$. Let $x=VS(b)$. Since
$\bar{C}_1C_i^{-1}\in G$ and $\bar{C}_1D_j^{-1}\in G$, we have 
$\bar{C}_1^{\, g-3}V\in G$. By the above lemmas, we also have $SBS^{-1}\bar{C}_1^{-1}\in G$. Therefore 
\begin{eqnarray} \label{eqn:vsb}
(\bar{C}_1^{ \, g-3}V )( SBS^{-1}\bar{C}_1^{-1})( V^{-1} \bar{C}_1^{\, 3-g})
&=&(V S)B( VS)^{-1} \bar{C}_1^{-1}\\
&=& X\bar{C}_1^{\, -1}
\end{eqnarray}
is contained in $G$. Therefore, $x\in \b$. Moreover, since $BD_2^{-1}\in G$ and $XD_2^{-1}\in G$,
we obtain 
$$
(BD_2^{-1})(XD_2^{-1})(D_2B^{-1})
= (BD_2^{-1})X (BD_2^{-1})^{-1}  D_2^{-1}
=A_4 D_2^{-1}
$$
is contained in $G$. This shows that $a_4$, and hence, all $a_i$, is in $\b$.

Suppose next that $g$ is odd and $g\geq 5$. Now again let $V$ denote the diffeomorphism
$
D_3^{-1}C_4^{-1}D_5^{-1}C_6^{-1}\cdots D_{g-2}^{-1}C_{g-1}^{-1}.
$
We have shown in the proof of Theorem~\ref{thmmcg} that the diffeomorphism $C_2D_1BVS$ maps the the curve $b$ 
to $a_3$. In this case, let $x$ denote the curve $VS(b)$, as above. The equation~(\ref{eqn:vsb}) shows that $x\in \b$.
We now use the equality
$$
(C_2D_1BC_4^{-3}) (XC_4^{-1}) (C_2D_1BC_4^{-3})^{-1}=A_3C_4^{-1}
$$
to conclude that $a_3$ is in $\b$.

Suppose now that $g=3$. It can easily be shown that 
$C_1^{-1} D_1^{-1} S^2 B S (b)=a_5$. We use this fact to prove that 
$a_5$ is contained in $\b$. Notice  that $C_2^{-1}S(b)=b_2$. The equation
$$
(\bar{C}_1 C_2^{\, -1}) ( SBS^{-1} C_2^{\, -1}) (C_2 \bar{C}_1^{\, -1})
=(C_2^{\, -1} S)B( S^{\, -1} C_2) C_2^{\, -1}
=B_2 C_2^{\, -1}
$$
shows that $b_2\in \b$. Let $y$ denote the curve $BS(b)$. From the equation 
$$
(B \bar{C}_1^{\, -1}) ( SBS^{-1} B_2^{\, -1}) (B \bar{C}_1^{\, -1})^{-1}=Y B_2^{\, -1},
$$
we conclude that $y\in \b$. By Lemma~\ref{lemma9}, $z=S^2(y)=S^2BS(b)$ is also contained in $\b$.
Finally, the fact that $a_5$ is in $\b$ follows from the equation
$$
(C_1^{\, -1}D_1^{\, -1}B_2^{\, 2}) ( Z  B_2^{\, -1})  
(B_2^{\, -2}D_1C_1)=A_5 B_2^{\, -1}.
$$

This concludes the proof of $a_i\in\b$. In order to finish the proof of the lemma,
it remains to prove that $b_2$ and $e$ are contained in $\b$. 

It is easy to see that $C_2^{-1}A_6A_5A_4(b)=b_2$ and $A_2A_1A_4^{-1}C_1(a_5)=e$.
It can be shown that the diffeomorphisms
$ A_1^{-2} C_2^{-1}A_6A_5A_4 $ and $B_2^{-2}A_2A_1A_4^{-1}C_1 $ are in $G$.
Finally, $b_2\in G$ follows from
$$
( A_1^{-2} C_2^{-1}A_6A_5A_4 ) (B A_1^{-1}) ( A_1^{-2} C_2^{-1}A_6A_5A_4)^{-1}
=B_2A_1^{-1},
$$
and  $e\in G$ follows from 
$$ 
( B_2^{-2}A_2A_1A_4^{-1}C_1 )
A_5B_2^{-1} ( B_2^{-2}A_2A_1A_4^{-1}C_1 )^{-1}
=EB_2^{-1}.
$$

This completes the proof.
\end{proof}

\begin{theorem}
The mapping class group $\mod_{g,1}$ $($and hence $\mod_g)$ is generated by two elements
of finite order. 
\end{theorem}
\begin{proof}
If $g=1$ then $A_2A_1$ and $A_1A_2A_1$ are of orders $6$ and $4$, respectively, and they
generate $\mod_{1,1}$. 
If $g=2$ then $A_4A_3A_2A_1$ and $A_5A_4A_3A_2A_1$ are of orders $10$ and $6$, respectively, and they
generate $\mod_{2,1}$. Suppose that $g\geq 3$ and let $G$ be the subgroup of $\mod_{g,1}$ generated by
the elements $S=A_{2g}A_{2g-1}\cdots A_2A_1$ and $BSB^{-1}$, which are both of order $4g+2$.
Since the curves $a_3,a_5,b, b_2, d_1$ and $e$ are all contained in $\b$, the elements
$BA_3^{-1}$, $ D_1A_5^{-1}$ and $ EB_2^{-1}$ are contained in $G$. By the lantern relation~(\ref{lant2}),
$A_1\in G$. Therefore, $S^{-i+1}A_1S^{i-1}=A_i \in G$.

Finally, the element $C_1^{-1}A_1^{-1}A_4A_3$ is in $G$ and maps $a_2$ to $b$. Since $A_2\in G$, this
shows that $B\in G$. Consequently, $G=\mod_{g,1}$.

This finishes the proof of the theorem.
\end{proof}

\section{A presentation of the mapping class group on two generators}

In this last section we transform the Wajnryb presentation of the mapping class group 
$\mod_g^1$ to a presentatoin on the two generators
$B$ and $S$. It turns out that the number of relations in the new presentation
depends linearly on $g$, whereas it is quadratic in the original presentation. 

\begin{theorem} $($\cite{w0,w2}$)$ \label{presentation} 
 Let $g\geq 3$.
 The mapping class group $\mod _g^1$ admits a presentation with generators
 $B,A_1,A_2,\ldots, A_{2g}$ and with defining relations
 \begin{eqnarray*}
 ({\rm i}) && BA_i=A_iB \ \ if \ \ i\neq 4,  \ and \ A_jA_k=A_kA_j
        \ if  \ |j-k|\geq 2, \\
 ({\rm ii}) && BA_4B=A_4BA_4, \  A_iA_{i+1}A_i=A_{i+1}A_iA_{i+1}
        \ for \ 1\leq i\leq 2g-1,\\
 ({\rm iii}) && (A_1A_2A_3)^4=B
        (A_4A_3A_2A_1^2A_2A_3A_4)B (A_4A_3A_2A_1^2A_2A_3A_4)^{-1},\\
 ({\rm iv}) && A_1 A_3 A_5 w B w^{-1} =
        (t_2t_1)^{-1} B (t_2t_1)  t_2^{-1} B t_2 B,
 \end{eqnarray*}
where
 \\
 $t_1=A_2 A_1 A_3 A_2, \,\, t_2=A_4 A_3 A_5 A_4$ and
 \\
 $w=A_6 A_5 A_4 A_3 A_2 (t_2 A_6 A_5)^{-1} B (t_2 A_6 A_5)
 (A_4A_3A_2A_1)^{-1}$. 
\end{theorem}

Suppose that $g\geq 3$.
We have shown above that the curve $b$ is mapped to $a_3$ by the diffeomorphism
\begin{eqnarray}\label{eqn:odd}
X_1=C_2D_1BD_3^{-1}C_4^{-1}D_5^{-1}C_6^{-1}\cdots D_{g-2}^{-1}C_{g-1}^{-1}S
\end{eqnarray}
if $g$ is odd, and to $a_4$ by
\begin{eqnarray}\label{eqn:even}
X_2=BD_2^{-1}C_3^{-1}D_4^{-1}C_5^{-1}\cdots D_{g-2}^{-1}C_{g-1}^{-1}S
\end{eqnarray}
if $g$ is even. 

From the equalities $S^{-2k+1}(b)=c_k$ and $S^{-2k}(b)=d_k$
we get that
$$S^{-2k+1}BS^{2k-1}=C_k$$ and
$$S^{-2k}BS^{2k}=D_k,$$
which should replace $C_k$ and $D_k$ in the
equations~(\ref{eqn:odd}) and~(\ref{eqn:even}). Also note that $S(a_k)=a_{k-1}$ 
and hence $SA_kS^{-1}=A_{k-1}$. 
Let us now define
\begin{eqnarray}\label{eqn:X}
X=\left\{
\begin{array}{ll}
S^{3}X_1 & \mbox{ if $g$ is odd,}\\
S^{4}X_2 & \mbox{ if $g$ is even,}
\end{array}
\right.
\end{eqnarray}
so that $S^{-k}X(b)=a_k$, and therefore 
\begin{eqnarray}\label{eqn:A_k}
A_k=S^{-k}XBX^{-1}S^{k}=B^{S^{-k}X}
\end{eqnarray}
for every $k=1,2,\ldots,2g$.

The presentation of the mapping class group $\mod_g^1$ in Theorem~\ref{presentation} 
has $2g+1$ generators and $2g^2+g+2$ relations.
By replacing each $A_k$ in Theorem~\ref{presentation} by $S^{-k}XBX^{-1}S^{k}$, we get
another presentation of $\mod_g^1$ on the two generators $B$ and $S$.
Now these relations are not as nice as those of Wajnryb's.
But, the number of relations in the new presentation
reduces to $4g+1$.

For any $i,j$ with $|i-j|>1$ and for any $k$, all the relations
$A_{i+k}A_{j+k}= A_{j+k}A_{i+k}$ reduce to the single relation
$B^X B^{S^{i-j}X} =B^{S^{i-j}X} B^X$. This shows that the relations 
in~(i) gives rise to $(2g-1)+(2g-2)=4g-3$ relations.

Each braid relation $A_kA_{k+1}A_k=A_{k+1}A_kA_{k+1}$ reduces to the single
relation $B^XB^{S^{-1}X}B^X =B^{S^{-1}X} B^X B^{S^{-1}X}$. Thus the relations in~(ii)
reduces to two relations.

Consequently the new presentation of $\mod_g^1$ is given as follows:
Let $P$ denote $B^X$ and let $U$ denote $S^{-1}$, where $X$ is defined as in~(\ref{eqn:X}).

\begin{theorem} \label{newpresentation} 
 Let $g\geq 3$.
 The mapping class group $\mod _g^1$ admits a presentation with generators
 $B$ and $U$,  and with defining relations
 
$({\rm i})$  \ \,$BP^{U^i}=P^{U^i}B, \mbox{ for }1\leq i\leq 2g \mbox{ and } i\neq 4$,

$({\rm ii})$  \,$PP^{U^i}=P^{U^i}P \mbox{ for }2\leq i\leq 2g-1$,

$({\rm iii})$   $BP^{U^4}B=P^{U^4}BP^{U^4} \mbox{ and } PP^UP=P^UPP^U,$

$({\rm iv})$   $(P^U P^{U^2}P^{U^3})^4=BB^V, \mbox{ and}$ 

$({\rm v})$   \,$P^UP^{U^3}P^{U^5}B^W=B B^{t_1^{-1}t_2^{-1}} B^{t_2^{-1}}$ , \\
where\\
$V=P^{U^4}P^{U^3}P^{U^2}P^{U^1}P^{U^1}P^{U^2}P^{U^3}P^{U^4}$.\\
$t_1= P^{U^2}P^{U}P^{U^3}P^{U^2}$, $t_2=P^{U^4}P^{U^3}P^{U^5}P^{U^4}$ and \\
$W=P^{U^6}P^{U^5}P^{U^4}P^{U^3}P^{U^2}(t_2P^{U^6}P^{U^5})^{-1}
B (t_2P^{U^6}P^{U^5}) (P^{U^4}P^{U^3}P^{U^2}P^U)^{-1}$.
\end{theorem}

The presentation of the mapping class group $\mod_g$ is obtained from that of
$\mod_g^1$ by adding one more relation (c.f. \cite{w2}). Hence, in the new presentation on the generators 
$B$ and $S$ the number of relation reduces to  $4g+2$.

\end{document}